\documentclass{article}
\usepackage{amsfonts}
\usepackage{amssymb}
\usepackage{graphicx}
\usepackage{amsmath}
\usepackage[singlespacing]{setspace}
\usepackage{fancyhdr}

\newtheorem{theorem}{Theorem}[section]

\newtheorem{definition}[theorem]{Definition}
\newtheorem{example}[theorem]{Example}

\newtheorem{lemma}[theorem]{Lemma}
\newtheorem{notation}[theorem]{Notation}

\newtheorem{proposition}[theorem]{Proposition}
\newtheorem{remark}[theorem]{Remark}

\newenvironment{proof}[1][Proof]{\noindent\textbf{#1.} 
}{\ \rule{0.6em}{0.6em}}
\input tcilatex 
\begin{document}

\title{A Minimal Groebner Basis for the Defining Ideals of Certain Affine
Monomial Curves}
\author{Ibrahim Al-Ayyoub}
\maketitle

\begin{abstract}
Let $K$ be a field and let $m_{0},...,m_{n}$ be an almost arithmetic
sequence of positive integers. Let $C$ \ be a monomial curve in the affine $%
\left( n+1\right) $-space, defined parametrically by $x_{0}=t^{m_{0}},\ldots
,x_{n}=t^{m_{n}}$. In this article we produce Groebner bases for the
defining ideal of $C$, correcting previous work of Sengupta,(2003).
\end{abstract}

\section*{Introduction}

Let $n\geq 2$, $K$ a field and let $x_{0},\ldots ,x_{n},t$ be
indeterminates. Let $m_{0},\ldots ,m_{n}$ be an almost arithmetic sequence
of positive integers, that is, some $n-1$ of these form an arithmetic
sequence, and assume $gcd(m_{0},\ldots ,m_{n})=1$. Let $P$\ be the kernel of
the $K$-algebra homomorphism $\eta :K[x_{0},\ldots ,x_{n}]\rightarrow K[t]$,
defined by $\eta (x_{i})=t^{m_{i}}$. A set of generators for the ideal $P$
was explicitly constructed in Patil and Singh (1990). We call these
generators the \textit{\textquotedblleft Patil-Singh
generators\textquotedblright }. Out of this generating set, Patil (1993)
constructed a minimal set $\Omega $ for the ideal $P$. We call the elements
of $\Omega $ \textit{\textquotedblleft Patil generators\textquotedblright }.
Sengupta (2003) proved that $\Omega $ forms a Groebner basis for the
relation ideal $P$ with respect to the grevlex monomial order, however,
Sengupta's proof is not complete, as in fact $\Omega $ is not a Groebner
basis in all cases, see Remark~\ref{Patil-Not-Gb} and Remark~\ref%
{PnotGBinSomeA3}. The goal of this article is to produce a minimal Groebner
basis for $P$. Remark~\ref{Patil-Not-Gb} was a motive to modify Patil
generators before computing Groebner basis. In Section~\ref{I&J} the set $%
\Omega $ is described more explicitly and it is modified. Then in Section~%
\ref{PtoGB} we state and prove a general result that we use to produce a
Groebner basis. This makes the proof of the main result of this thesis much
shorter and simpler than the work given in Sengupta (2003).

\section{Patil Generators\label{Patil-Gens}}

We shall use the notation and the terminology from Patil and Singh (1990)
and Patil (1993) with a slight difference in naming some variables and
constants. Let $n\geq 2$ be an integer and let $p=n-1$ . Let $m_{0},\ldots
,m_{p}$ be an arithmetic sequence of positive integers with $0<m_{0}<\cdots
<m_{p}$, let $m_{n}$ be arbitrary, and $gcd(m_{0},\ldots ,m_{n})=1$. Let $%
\Gamma $ denote the numerical semigroup that is generated by $m_{0},\ldots
,m_{n}$ i.e. $\Gamma =\sum\limits_{i=0}^{n}\mathbb{N}_{\mathbf{0}}m_{i}$. We
assume throughout that $\Gamma $\ is minimally generated by $m_{0},\ldots
,m_{n}$. Put $\Gamma ^{\prime }=\sum\limits_{i=0}^{p}\mathbb{N}_{0}m_{i}$.
Thus $\Gamma =\Gamma ^{\prime }+\mathbb{N}_{\mathbf{0}}m_{n}$.

\begin{notation}
\label{g-sub-t}\textit{For }$a,b\in \mathbb{Z}~$\textit{\ let }$[a,b]=\{t\in 
\mathbb{Z}\mid a\leq t\leq b\}$\textit{. For }$t\geq 0$\textit{, let }$%
q_{t}\in \mathbb{Z}\mathbf{,}$\textit{\ }$r_{t}\in \lbrack 1,p]$\textit{\
and }$g_{t}\in \Gamma ^{\prime }$\textit{\ \ be defined by }$t=q_{t}p+r_{t}$%
\textit{\ and }$g_{t}=q_{t}m_{p}+m_{r_{t}}$.\ \ \ \ \ \ \ \ \ \ \ \ \ 
\end{notation}

Let $S=\{\gamma \in \Gamma \mid \gamma -m_{0}\notin \Gamma \}$. As all large
integers are in $\Gamma $, $S$ is a finite set of non-negative integers. The
following gives an explicit description of $S$.\ \ \ \ \ 

\begin{lemma}
\label{Parameters}(Lemma 3.1, 3.2 in Patil-Singh (1990))\textbf{\ }\textit{%
Let }$u=min\{t\geq 0\mid g_{t}\notin S\}$\textit{\ and }$\upsilon
=min\{b\geq 1\mid bm_{n}\in \Gamma ^{\prime }\}$\textit{.}\newline
(a) \textit{There exist unique integers }$w\in \lbrack 0,\upsilon -1]$%
\textit{, }$z\in \lbrack 0,u-1]$\textit{, }$\lambda \geq 1$\textit{, }$\mu
\geq 0$\textit{, and }$\nu \geq 2$\textit{\ such that\newline
\ \ (i) }$g_{u}=\lambda m_{0}+wm_{n}$\textit{;\newline
\ \ (ii) }$\upsilon m_{n}=\mu m_{0}+g_{z}$\textit{;\newline
\ \ (iii) }$g_{u-z}+(\upsilon -w)m_{n}=\nu m_{0}$\textit{, where }$\nu
=\left\{ 
\begin{tabular}{ll}
$\lambda +\mu +1\text{,}$ & if$\text{\ \ }r_{u-z}<r_{u}\text{;}$ \\ 
$\lambda +\mu \text{,}$ & $\text{if \ }r_{u-z}\geq r_{u}\text{.}$%
\end{tabular}%
\right. $\newline
\newline
\textit{(b) Let }$\ V=[0,u-1]\times \lbrack 0,\upsilon -1]$\textit{\ and }$%
W=[u-z,u-1]\times \lbrack \upsilon -w,\upsilon -1]$\textit{. Then every
element of }$\Gamma $ can be expressed uniquely in the form $%
am_{0}+g_{s}+bm_{n}$ with $a\in \mathbb{N}_{0}$ and $(s,b)\in V-W.$
\end{lemma}

\begin{proof}
The proof can be found in Patil and Singh (1990). We would like to recall
their proof for part (b) which constitutes of a major point in our proof of
the main theorem of this article. The existence of $a\in \mathbb{N}_{0}$ and 
$(s,b)\in V-W$ such that every element of $\Gamma $ can be expressed in the
form $am_{0}+g_{s}+bm_{n}$ is a consequence of the main proposition of this
article, Proposition~\ref{MainProp}. Here we prove uniqueness. Let $\gamma
\in \Gamma $ such that $\gamma =am_{0}+g_{s}+bm_{n}$ and also $\gamma
=a^{\prime }m_{0}+g_{t}+cm_{n}$ with $(s,b),(t,c)\in V-W.$ This implies $%
g_{s}+bm_{n}\equiv g_{t}+cm_{n}$ $\left( \text{mod }m_{0}\right) $. Assume $%
b\geq c$, then $0\leq (b-c)m_{n}\equiv g_{s-t}$ $\left( \text{mod }%
m_{0}\right) $, hence we also must assume that $s\leq t$. We consider two
cases: let $e\geq 0.$ Case $(b-c)m_{n}=g_{s-t}+em_{0}$ is a contradiction
since $b-c<\upsilon $. Case $(b-c)m_{n}+em_{0}=g_{s-t}$. If $e>0$ then $%
g_{s-t}-m_{0}=(b-c)m_{n}+(e-1)m_{0}\in \Gamma $, but this is a contradiction
since $s-t<u$. Thus, $e=0$. Therefore, $(b-c)m_{n}=g_{s-t}$. But $%
b-c<\upsilon $ and $s-t<u$ by assumption. Hence, $b=c$ and $s=t$.
\end{proof}

\begin{notation}
\label{W&q'&r'}\textit{Let }$q=q_{u},$\textit{\ }$r=r_{u},$\textit{\ }$%
q^{\prime }=q_{u-z},$\textit{\ }$r^{\prime }=r_{u-z}$\textit{. From now on,
the symbols }$q,$\textit{\ }$q\prime ,$\textit{\ }$r,$\textit{\ }$r\prime ,$%
\textit{\ }$u,\upsilon ,$\textit{\ }$w,$\textit{\ }$z,$\textit{\ }$\lambda ,$%
\textit{\ }$\mu ,$\textit{\ }$\nu ,V$\textit{\ and }$W$\textit{\ will have
the meaning assigned to them by this notation and the lemma above.}
\end{notation}

\begin{remark}
\label{q>0-and-u>p}Note that for $1\leq i\leq p$ we have $%
g_{i}-m_{0}=m_{i}-m_{0}.$ \textit{Then by the minimality assumption on the
generators of \ }$\Gamma $\textit{\ it follows that} $u>p$, hence $q>0.$
\end{remark}

We recall the construction and the result given in Patil (1993): recall $%
p=n-1$ and let

$\xi _{i,j}=\left\{ 
\begin{tabular}{ll}
$x_{i}x_{j}-x_{0}x_{i+j}$, & $\text{if\ \ \ \ }i+j\leq p$; \\ 
$x_{i}x_{j}-x_{i+j-p}x_{p}$, & $\text{if \ \ \ }i+j>p$,%
\end{tabular}%
\right. $

$\varphi _{i}=x_{r+i}x_{p}^{q}-x_{0}^{\lambda -1}x_{i}x_{n}^{w}$,

$\psi _{j}=x_{r^{\prime }+j}x_{p}^{q^{\prime }}x_{n}^{\upsilon
-w}-x_{0}^{\nu -1}x_{j}$,

$\theta $\ $=\left\{ 
\begin{tabular}{ll}
$x_{n}^{\upsilon }-x_{0}^{\mu }x_{r-r^{\prime }}x_{p}^{q-q^{\prime }}\text{,}
$ & if $\text{\ }r^{\prime }<r\text{;}$ \\ 
$x_{n}^{\upsilon }-x_{0}^{\mu }x_{p+r-r^{\prime }}x_{p}^{q-q^{\prime }-1}%
\text{,}$ & if $\text{\ }r^{\prime }\geq r$,%
\end{tabular}%
\right. $\newline
Let

$I=\left\{ 
\begin{tabular}{ll}
$\lbrack 0,p-r]\text{,}$ & if$\text{ }\mu \neq 0\text{ }$or$\text{ }W=\phi 
\text{;}$ \\ 
$\lbrack \max (r_{z}-r+1,0),p-r]\text{,}$ & if$\text{ }\mu =0\text{ }$and$%
\text{ }W\neq \phi $,%
\end{tabular}%
\ \right. $

$J=\left\{ 
\begin{tabular}{ll}
$\phi \text{,}$ & if $\text{\ }W=\phi \text{;}$ \\ 
$\lbrack 0,\min (z-1,p-r^{\prime })]\text{,}$ & if$\text{ \ }W\neq \phi 
\text{.}$%
\end{tabular}%
\ \right. $

\begin{theorem}
\label{Patil&Patil-Singh Gens}(Theorem (4.5) in Patil (1993), and Theorem
(4.5) in Patil-Singh (1990)) The set%
\begin{equation*}
\Omega =\{\xi _{i,j}\mid 1\leq i\leq j\leq p-1\}\cup \{\theta \}\cup
\{\varphi _{i}\mid i\in I\}\cup \{\psi _{j}\mid j\in J\}
\end{equation*}%
forms a minimal generating set for the ideal $P$. Also, the set%
\begin{equation*}
\{\xi _{i,j}\mid 1\leq i\leq j\leq p-1\}\cup \{\theta \}\cup \{\varphi
_{i}\mid 0\leq i\leq p-r\}\cup \{\psi _{j}\mid 0\leq j\leq p-r^{\prime }\}
\end{equation*}%
forms a generating set for the ideal $P$.
\end{theorem}

Sengupta (2003) tried to prove that the set $\Omega $ from Theorem~\ref%
{Patil&Patil-Singh Gens} forms a Groebner basis for $P$ with respect to the
grevlex monomial order using the grading $wt(x_{i})=m_{i}$ with\textit{\ }$%
x_{0}<x_{1}<\cdots <x_{n}$. In this case \textit{\ }$\prod%
\limits_{i=0}^{n}x_{i}^{a_{i}}>_{grevlex}\prod\limits_{i=0}^{n}x_{i}^{b_{i}}$%
\ if in the ordered tuple $(a_{1}-b_{1},\ldots ,a_{n}-b_{n})$\ the left-most
nonzero entry is negative. Sengupta's proof works for arithmetic sequences,
but it is incomplete for the almost arithmetic sequences:\ \ \ 

\begin{remark}
\label{Patil-Not-Gb}Assume $r^{\prime }\geq r$ , $\mu =0$, and$\text{ }W\neq
\phi $. Then \textit{Patil generators are not a Groebner basis with respect
to the grevlex monomial ordering with }$x_{0}<x_{1}<\cdots <x_{n}$ \textit{%
and\ with the grading }$wt(x_{i})=m_{i}$\textit{.}
\end{remark}

\begin{proof}
As $u-z=(q-q_{z})p+(r-r_{z})$ then $r^{\prime }\geq r$ if and only if $%
r_{z}\geq r$. Assume $r^{\prime }\geq r$, then $r_{z}-r+1>0$ and also $%
\theta =x_{n}^{\upsilon }-x_{0}^{\mu }x_{p+r-r^{\prime }}x_{p}^{q-q^{\prime
}-1}$. Assume also that $\mu =0\text{ }$and$\text{ }W\neq \phi $, then $%
I=[\max (r_{z}-r+1,0),p-r]=[r_{z}-r+1,p-r]$. Under these assumptions the
S-polynomial $S(\psi _{k},\theta )$\ can not be reduced to zero modulo $%
\Omega $: for $0\leq k<r_{z}-r+1$ consider $S(\psi _{k},\theta )=x_{0}^{\mu
}S_{1}$ where $S_{1}=x_{0}^{\lambda -1}x_{k}x_{n}^{w}-\underline{%
x_{r^{\prime }+k}x_{p+r-r^{\prime }}x_{p}^{q-1}},$ with the leading monomial
underlined. We note that $LM(S_{1})$ is a multiple of $LM(\xi _{r^{\prime
}+j,p+r-r^{\prime }})$ only. Hence, the only possible way to reduce $S_{1}$
with respect to $\Omega $ is by using $\xi _{r^{\prime }+j,p+r-r^{\prime }}$%
. However, none of the terms of the binomial $S_{1}+x_{p}^{q-1}\xi
_{r^{\prime }+j,p+r-r^{\prime }}=x_{r+k}x_{p}^{q}-x_{0}^{\lambda
-1}x_{k}x_{n}^{w}$ is a multiple of any of the leading terms of Patil
generators. Therefore, it can not be reduced to $0$ modulo $\Omega $.
\end{proof}

It remains to give a concrete example for which $r^{\prime }\geq r$ , $\mu
=0 $, and$\text{ }W\neq \phi $ so that Patil generators are not a Groebner
basis:

\begin{remark}
\label{PnotGBinSomeA3}Let $m_{0}\geq 5$ be an odd integer. Let $P$ be the
defining ideal of the monomial curve that corresponds to the almost
arithmetic sequence $m_{0},m_{0}+1,m_{0}-1$. Then Patil generators for the
ideal $P$ are not\ a Groebner basis \textit{with respect to the grevlex
monomial ordering with }$x_{0}<x_{1}<\cdots <x_{n}$ \textit{and with the
grading }$wt(x_{i})=m_{i}$\textit{.}
\end{remark}

\begin{proof}
Observe: $p=1,n=2$ , and $g_{i}=i(m_{0}+1)$ for all $i$.

Let $\upsilon ,\mu ,$ and $z$\ be as defined in Lemma~\ref{Parameters}. Then 
$\upsilon (m_{0}-1)=\mu m_{0}+z(m_{0}+1)$ for some integers $\mu ,z\geq 0$ .
This implies $\mu +z<v.$ Note that%
\begin{equation*}
\upsilon (m_{0}-1)=\mu m_{0}+z(m_{0}+1)=(\mu +z)(m_{0}-1)+\mu +2z.
\end{equation*}%
Thus $\mu +2z=s(m_{0}-1)$ for some $s\geq 1$. Hence,%
\begin{equation*}
\upsilon >\mu +z\geq \dfrac{\mu }{2}+z=\dfrac{s}{2}(m_{0}-1)\geq \frac{%
m_{0}-1}{2}.
\end{equation*}%
Thus,%
\begin{equation}
\upsilon \geq \frac{m_{0}+1}{2}.  \label{v>or=(m0+1)/2}
\end{equation}%
On the other hand, note that%
\begin{equation}
\dfrac{m_{0}+1}{2}(m_{0}-1)=\dfrac{m_{0}-1}{2}(m_{0}+1)\in \Gamma ^{\prime }.
\label{v((m0-1))}
\end{equation}%
Therefore, by the minimality of $\upsilon $ we must have%
\begin{equation}
\upsilon \leq \frac{m_{0}+1}{2}.  \label{v<=(m0+1)/2}
\end{equation}%
By~$\left( \ref{v>or=(m0+1)/2}\right) $ and ~$\left( \ref{v<=(m0+1)/2}%
\right) $ it follows that $\upsilon =\dfrac{m_{0}+1}{2}$.

Let $u,\lambda ,w,$ and $g_{u}$ be as defined in Lemma~\ref{Parameters}. Note%
\begin{equation}
\dfrac{m_{0}+1}{2}(m_{0}+1)-m_{0}=\dfrac{m_{0}-1}{2}(m_{0}-1)+m_{0}\in
\Gamma .  \label{u((m0+1))}
\end{equation}%
Therefore,%
\begin{equation}
u\leq \dfrac{m_{0}+1}{2}.  \label{u<or=(m0+1)/2}
\end{equation}

Claim $w>0$: if $w=0$ then $g_{u}=\lambda m_{0}$, thus $u(m_{0}+1)=\lambda
m_{0}$. But $m_{0}$ and $m_{0}+1$ are relatively prime, therefore, we must
have $u=bm_{0}$ for some $b\geq 1$, a contradiction to~$\left( \ref%
{u<or=(m0+1)/2}\right) $. Thus $w>0$.

Claim $\lambda <u$: by Lemma~\ref{Parameters} we have $u(m_{0}+1)=\lambda
m_{0}+w(m_{0}-1)$. If $\lambda \geq u$ then $w(m_{0}-1)=u(m_{0}+1)-\lambda
m_{0}=u+(u-\lambda )m_{0}$, which implies $u\geq m_{0}-1$ as $w>0$, a
contradiction to~$\left( \ref{u<or=(m0+1)/2}\right) $. Thus Claim $\lambda
<u $.

Now consider%
\begin{eqnarray*}
w(m_{0}-1) &=&u(m_{0}+1)-\lambda m_{0}, \\
&=&\left( u-\lambda \right) (m_{0}-1)+2u-\lambda .
\end{eqnarray*}%
As $w(m_{0}-1)>0$ and $u>\lambda $ we must have $2u-\lambda =c(m_{0}-1)$ for
some $c\geq 1$. But if $u\leq \dfrac{m_{0}-1}{2}$ then $2u-\lambda \leq
m_{0}-1-\lambda $, a contradiction as $\lambda \geq 1$. Therefore, 
\begin{equation}
u>\dfrac{m_{0}-1}{2}.  \label{u>>(m0-1)/2}
\end{equation}%
By $\left( \ref{u<or=(m0+1)/2}\right) $ and $\left( \ref{u>>(m0-1)/2}\right) 
$ it follows that $u=\dfrac{m_{0}+1}{2}$.

Now by the uniqueness in Lemma~\ref{Parameters} and as of $\left( \ref%
{v((m0-1))}\right) $ and $\left( \ref{u((m0+1))}\right) $ it follows that $%
\mu =0$, $z=\dfrac{m_{0}-1}{2}$, $\lambda =2$ and $w=\dfrac{m_{0}-1}{2}$.

Finally, note that $r=p=r^{\prime }=1$. Therefore, the parameters $z,w,\mu
,p,r,$ and $r^{\prime }$ all satisfy the assumptions of the previous remark,
hence done.
\end{proof}

\begin{example}
\label{7,8,6 not GB}Let $m_{0}=7,m_{1}=8,$ and $m_{2}=6$. By the previous
remark $\upsilon =4$, $z=3$, $\mu =0$, $u=4$, $\lambda =2$, and $w=3$.
Therefore, as $p=1$ then $q=3,r=1,q^{\prime }=0,r^{\prime }=1,q_{z}=2$, and $%
r_{z}=1$. Thus $I=[\max (r_{z}-r+1,0),p-r]=\phi $ and $\ J=[0,\min
(z-1,p-r^{\prime })]=\{0\}$. By \ Patil (1993), or Theorem~\ref{Patil-Gens},
a minimal set of generators for the ideal $P$\ is $\Omega =\{\psi
_{0},\theta \}$\ where $\psi _{0}=\underline{x_{1}x_{2}}-x_{0}^{2}$ and $%
\theta =\underline{x_{2}^{4}}-x_{1}^{3}$. This is not a Groebner basis with
respect the grevlex monomial order with $x_{0}<x_{1}<x_{2}$ and with the
grading $wt(x_{i})=m_{i}$: note neither term of $S(\psi _{0},\theta
)=x_{1}^{4}-x_{0}^{2}x_{2}^{3}$ is a multiple of $LT(\psi _{0})\ $or $%
LT(\theta )$.
\end{example}

\section{Explicit Descriptions for I and J\label{I&J}}

Let $u,$ $z,$ $q,$ $r,$ $q^{\prime }=q_{u-z}$, and $r^{\prime }=r_{u-z}$ be
as in Lemma~\ref{Parameters} and Notation~\ref{W&q'&r'}. Let $z=q_{z}p+r_{z}$
with $q_{z}\in \mathbb{Z}$ and $r_{z}\in \lbrack 1,p]$. By Notation~\ref%
{g-sub-t} it is clear that $q_{z}\leq q$ since $0\leq z\leq u-1$. As $%
u-z=(q-q_{z})p+(r-r_{z})$, it follows that $q^{\prime }=q-q_{z}-\varepsilon $
and $r^{\prime }=\varepsilon p+r-r_{z}$ where $\varepsilon =\left\{ 
\begin{tabular}{ll}
$0\text{,}$ & if$\text{ \ }r>r_{z}\text{;}$ \\ 
$1\text{,}$ & if$\text{ \ }r\leq r_{z}$.%
\end{tabular}%
\right. $

Therefore, $r^{\prime }<r$\ if and only if $r_{z}<r$. We note the following:
if $r^{\prime }<r$ then $r>r_{z}$, hence $r-r^{\prime }=r_{z}$ and $%
q-q^{\prime }=q_{z}$. Also, if $r^{\prime }\geq r$ then $r\leq r_{z}$, hence 
$p+r-r^{\prime }=r_{z}$ and $q-q^{\prime }-1=q_{z}$. Therefore,%
\begin{equation*}
\theta =x_{n}^{\upsilon }-x_{0}^{\mu }x_{r_{z}}x_{p}^{q_{z}}.
\end{equation*}

\begin{proposition}
\label{Min_z-1&p-r'}\textit{Let }$z>0$\textit{\ and let }$z=q_{z}p+r_{z%
\mathit{\ }}$\textit{with }$q_{z}\in \mathbb{Z\ }$\textit{and }$r_{z}\in
\lbrack 1,p]$\textit{. Then }$min\{z-1,p-r^{\prime }\}=\left\{ 
\begin{tabular}{ll}
$p-r^{\prime }\text{,}$ & if$\text{ \ }r\leq r_{z}\text{;}$ \\ 
$p-r^{\prime }\text{,}$ & if$\text{ \ }r>r_{z}\text{ and\ }z>p\text{;}$ \\ 
$z-1\text{,}$ & $\text{if \ }r>r_{z}\text{ and }z\leq p\text{.}$%
\end{tabular}%
\right. $\newline
\textit{Moreover, }$z\leq p$\textit{\ if and only if }$q_{z}=0$\textit{.}
\end{proposition}

\begin{proof}
First note that $p-r^{\prime }=(1-\varepsilon )p+r_{z}-r$. It is obvious
that if $z>0$\ then $q_{z}\geq 0$. Consider three cases:\newline
Case $r\leq r_{z}$: since $r\in \lbrack 1,p]$\ then $z-1=q_{z}p+r_{z}-1\geq
r_{z}-1\geq r_{z}-r=p-r^{\prime }$.\newline
Case $r>r_{z}$\ and $z>p$: this implies $q_{z}\geq 1.$\ Therefore, $%
z-1=q_{z}p+r_{z}-1\geq p+r_{z}-1\geq p+r_{z}-r=p-r^{\prime }$.\newline
Case $r>r_{z}$\ and $z\leq p$: this implies $q_{z}=0$. Therefore, $%
z-1=r_{z}-1\leq r_{z}-1+p-r<p+r_{z}-r=p-r^{\prime }$.
\end{proof}

Therefore, for $W\neq \phi $ we write $J$\ as following

$J=\left\{ 
\begin{tabular}{ll}
$\lbrack 0,p-r^{\prime }]\text{,}$ & if $q_{z}>0$ or $\varepsilon >0$; \\ 
$\lbrack 0,r_{z}-1]\text{,}$ & if $q_{z}=0$ and $\varepsilon =0$.%
\end{tabular}%
\right. $

\section{Groebner Basis\label{PtoGB}\ }

We will prove that the binomials listed here form a Groebner basis for the
defining ideal $P$ with respect to the grevlex monomial order with $%
x_{0}<x_{1}<\cdots <x_{n}$ and with the grading $wt(x_{i})=m_{i}$:

\ \ \ 

\ \ 
\begin{tabular}{lll}
$\varphi _{i}$ & $=\underline{x_{r+i}x_{p}^{q}}-x_{0}^{\lambda
-1}x_{i}x_{n}^{w}$, & for \ \ $\ 0\leq i\leq p-r$; \ \  \\ 
$\psi _{j}$ & $=\underline{x_{r^{\prime }+j}x_{p}^{q^{\prime
}}x_{n}^{\upsilon -w}}-x_{0}^{\nu -1}x_{j}$, & for \textit{\ \ \ }$j\in J$
\\ 
$\theta $ & $=\underline{x_{n}^{\upsilon }}-x_{0}^{\mu
}x_{r_{z}}x_{p}^{q_{z}}$, &  \\ 
$\xi _{i,j}$ & $=\left\{ 
\begin{tabular}{ll}
$\underline{x_{i}x_{j}}-x_{0}x_{i+j}$, & $\text{if\ \ \ \ }i+j\leq p$; \\ 
$\underline{x_{i}x_{j}}-x_{i+j-p}x_{p}$, & $\text{if \ \ \ }i+j>p$,%
\end{tabular}%
\right. $ & for\ \ \ $1\leq i\leq j\leq p-1$,%
\end{tabular}%
\newline
\ \ \newline
where the underlined monomials are the leading monomials.

Before we state the theorem that contains the main result of this article,
let us recall the definition of the minimal Groebner basis.\ \ \ \ \ \ \ \ \
\ \ \ \ \ \ \ \ \ \ 

\begin{definition}
\label{minimalGB}\textit{A \textbf{minimal Groebner basis} for a polynomial
ideal }$I$\textit{\ is a Groebner basis }$G$\textit{\ for }$I$\textit{\ such
that: \newline
(i) }$LC(f)=1$\textit{\ for all }$f\in G$\textit{, where }$LC(f)$\textit{\
is the leading coefficient of }$f$\textit{.\newline
(ii) For all }$f\in G$\textit{, }$LM(f)\notin \langle LM\{G-\{f\}\}\rangle ,$%
\textit{\ where }$LM\{G-\{f\}\}$ is the set of leading monomials of all
polynomials in the set $\{G-\{f\}\}$.\ \ \ \ \ 
\end{definition}

\ \ \ Now we state the main theorem of this article considering the
binomials $\varphi _{i}$, $\psi _{j}$, $\theta $, and $\xi _{i,j}$ as stated
above.

\begin{theorem}
\label{MainThm}\textit{The set }%
\begin{equation*}
\Phi =\mathit{\ }\{\varphi _{i}\mid 0\leq i\leq p-r\}\cup \{\psi _{j}\mid
j\in J\}\cup \{\theta \}\mathit{\ }\cup \mathit{\ }\{\xi _{i,j}\mid 1\leq
i\leq j\leq p-1\}
\end{equation*}%
\textit{is a minimal Groebner basis for the ideal }$P$\textit{\ with respect
to the grevlex monomial order with }$x_{0}<x_{1}<\cdots <x_{n}$\textit{\ and
with the grading }$wt(x_{i})=m_{i}$\textit{.}
\end{theorem}

The remainder of this article is to prove Theorem~\ref{MainThm}. The proof
is a consequence of a general result that we state and prove in the main
proposition of this article, Proposition~\ref{MainProp}. The following three
lemmas are needed to prove the main proposition.

\begin{lemma}
\label{Lemma1}Let $e_{0},\ldots ,e_{p},d_{1}\in \mathbb{N}_{0}$ with $%
e_{1}+\cdots +e_{p}\geq 1$. \textit{The monomial }$\alpha
=x_{0}^{e_{0}}\cdots x_{p}^{e_{p}}x_{n}^{d_{1}}$\textit{\ can }be reduced to 
$x_{0}^{h}x_{s}x_{p}^{l}x_{n}^{d}$\textit{\ (modulo }$\Phi $) \textit{with
respect to the grevlex monomial order }such that $l\leq q$, $d<\upsilon $, $%
s\in \lbrack 1,p]$, $h\geq e_{0}$, and with $\sum%
\limits_{i=0}^{p}e_{i}m_{i}+d_{1}m_{n}=hm_{0}+m_{s}+lm_{p}+dm_{n}$\textit{.}
\end{lemma}

\begin{proof}
Note that, using $\theta $ as many times as necessary, the monomial $\alpha
=x_{0}^{e_{0}}\cdots x_{p}^{e_{p}}x_{n}^{d_{1}}$ can be reduced to $%
x_{0}^{e_{0}^{\prime }}\cdots x_{p}^{e_{p}^{\prime }}x_{n}^{k}$\ with $%
k<\upsilon $ and $e_{1}^{\prime }+\cdots +e_{p}^{\prime }\geq e_{1}+\cdots
+e_{p}$. This reduces the proof to the case $\alpha =x_{0}^{e_{0}}\cdots
x_{p}^{e_{p}}x_{n}^{d_{1}}$ with $d_{1}<\upsilon $.

Using various $\xi _{i,j}$, $\alpha $\ can be reduced to $\sigma
=x_{0}^{h_{1}}x_{s_{1}}x_{p}^{l_{1}}x_{n}^{d_{1}}$ for some $h_{1}$, $%
l_{1}\in \mathbb{N}_{0}$, and $s_{1}\in \lbrack 1,p]$ with $%
\sum\limits_{i=0}^{p}e_{i}m_{i}=h_{1}m_{0}+m_{s_{1}}+l_{1}m_{p}$. Observe
that $g_{l_{1}p+s_{1}}=l_{1}m_{p}+m_{s_{1}}$.

If $l_{1}\leq q$ then there is nothing to prove. If $l_{1}>q$ apply the
following algorithm:\ 

\ \ \newline
Input $\sigma _{1}=x_{0}^{h_{1}}x_{s_{1}}x_{p}^{l_{1}}x_{n}^{d_{1}}$ .%
\newline
Output $\sigma =x_{0}^{h}x_{s}x_{p}^{l}x_{n}^{d}$ with $l\leq q$, $%
d<\upsilon $, and $s\in \lbrack 1,p]$ such that%
\begin{equation*}
hm_{0}+m_{s}+lm_{p}+dm_{n}=h_{1}m_{0}+m_{s_{1}}+l_{1}m_{p}+d_{1}m_{n}.
\end{equation*}%
Let $i=0$.\newline
\textbf{REPEAT }\newline
Let $i=i+1$.\newline
\textbf{Step 1}: Write $l_{i}=a_{i}(q+1)+b_{i}$ with $a_{i}\in \mathbb{N}%
_{0} $ and $b_{i}\in \lbrack 0,q].$\newline
Reduce $\sigma _{1}$ with respect to $\varphi _{p-r}$ $\ a_{i}$ times:%
\begin{equation*}
\sigma ^{\prime }=x_{0}^{a_{i}(\lambda
-1)+h_{i}}x_{s_{i}}x_{p-r}^{a_{i}}x_{p}^{b_{i}}x_{n}^{a_{i}w+d_{i}}.
\end{equation*}%
\textbf{Step 2}: Write $a_{i}w+d_{i}=c_{i}\upsilon +d_{i+1}$ with $c_{i}\in 
\mathbb{N}_{0}$ and $d_{i+1}\in \lbrack 0,\upsilon -1]$ (note that $%
c_{i}\leq a_{i}$ because $w<\upsilon $ and $d_{i}<\upsilon $. Also, note
that if $c_{i}=a_{i}$ then $d_{i+1}<d_{i}$). \newline
Reduce $\sigma ^{\prime }$ with respect to $\theta $ $\ c_{i}$ times: 
\begin{equation*}
\sigma ^{\prime \prime }=x_{0}^{a_{i}(\lambda -1)+h_{i}+c_{i}\mu
}x_{s_{i}}x_{r_{z}}^{c_{i}}x_{p-r}^{a_{i}}x_{p}^{c_{i}q_{z}+b_{i}}x_{n}^{d_{i+1}}.
\end{equation*}%
\textbf{Step 3}: Reduce $\sigma ^{\prime \prime }$ with respect to $\xi
_{r_{z},p-r}$ (or $\xi _{p-r,r_{z}}$ if $p-r\leq r_{z}$. Recall that $\xi
_{p-r,r_{z}}=x_{0}^{1-\varepsilon }x_{p-r^{\prime }}x_{p}^{\varepsilon }$
and $p-r^{\prime }=(1-\varepsilon )p+r_{z}-r$) $c_{i}$ times (recall $%
c_{i}\leq a_{i}$):%
\begin{eqnarray*}
\sigma ^{\prime \prime \prime } &=&x_{0}^{a_{i}(\lambda -1)+h_{i}+c_{i}\mu
}x_{s_{i}}\left( x_{0}^{1-\varepsilon }x_{p-r^{\prime }}x_{p}^{\varepsilon
}\right)
^{c_{i}}x_{p-r}^{a_{i}-c_{i}}x_{p}^{c_{i}q_{z}+b_{i}}x_{n}^{d_{i+1}}, \\
&=&x_{0}^{a_{i}(\lambda -1)+h_{i}+c_{i}\mu }x_{s_{i}}\left(
x_{0}^{1-\varepsilon }x_{p}^{\varepsilon }\right) ^{c_{i}}x_{p-r^{\prime
}}^{c_{i}}x_{p-r}^{a_{i}-c_{i}}x_{p}^{c_{i}q_{z}+b_{i}}x_{n}^{d_{i+1}}.
\end{eqnarray*}%
Reduce $\sigma ^{\prime \prime \prime }$ with respect to various $\xi _{i,j}$
we get:%
\begin{equation*}
\sigma _{i+1}=x_{0}^{h_{i+1}}x_{s_{i+1}}x_{p}^{l_{i+1}}x_{n}^{d_{i+1}}.
\end{equation*}%
with (by degree count) $l_{i+1}\leq c_{i}q_{z}+b_{i}+a_{i}+\varepsilon c_{i}$%
, and $h_{i+1}\geq a_{i}(\lambda -1)+h_{i}+c_{i}\mu $, and some $s_{i+1}\in
\lbrack 1,p].$\newline
\textbf{UNTIL}\ $l_{i+1}\leq q$.

\ \ \ 

Note that $l_{i+1}\leq l_{i}$ since $c_{i}\leq a_{i}$ and since%
\begin{equation}
l_{i+1}\leq c_{i}q_{z}+b_{i}+a_{i}+\varepsilon c_{i}\leq a_{i}(q-\varepsilon
)+b_{i}+(1+\varepsilon )a_{i}=a_{i}(q+1)+b_{i}=l_{i}.  \label{l_i+1=<l_i}
\end{equation}%
Now we prove that the sequence $l_{1},l_{2},\ldots $ is a decreasing
sequence which eventually goes below $q$:\ \newline
Case $c_{i}<a_{i}$, $q_{z}=0$ and $\varepsilon =0$: then%
\begin{equation*}
l_{i+1}\leq c_{i}q_{z}+b_{i}+a_{i}+\varepsilon
c_{i}=b_{i}+a_{i}<a_{i}(q+1)+b_{i}=l_{i}.
\end{equation*}%
Case $c_{i}<a_{i}$, $q_{z}\neq 0$ or $\varepsilon \neq 0$: then (recall that 
$q>0$ by Remark~\ref{q>0-and-u>p}) then $l_{i+1}\leq l_{i}$ by (\ref%
{l_i+1=<l_i}).\newline
Case $c_{i}=a_{i}$ for some $i>0$: then we prove that there must be an
integer $t>i$ such that $c_{t}<a_{t}$ and hence $l_{t+1}<l_{t}$: note that
if $c_{i}=a_{i}$ then $d_{i}>d_{i+1\newline
}$. Hence $d_{i},d_{i+1\newline
},\ldots $ is a strictly decreasing sequence of integers whenever $%
c_{i}=a_{i},c_{i+1}=a_{i+1},\ldots $ since the $d_{i}$ are non-negative hen
necessarily there exist $t>i$ such that $c_{t}<a_{t}$. Therefore, $%
l_{t+1}<l_{t}$.
\end{proof}

\begin{lemma}
\label{Lemma2}Assume $W\neq \phi $. Let $\sigma
=x_{0}^{h_{1}}x_{s_{1}}x_{p}^{l_{1}}x_{n}^{d_{1}}$ with $l_{1}>q^{\prime }$, 
$d_{1}\geq \upsilon -w$, $h_{1}\geq 0$, and $s_{1}\in \lbrack 1,p].$ Then $%
\sigma $ can be reduced \textit{(modulo }$\Phi $) \textit{with respect to
the grevlex monomial order }to $x_{0}^{h}x_{s}x_{p}^{l}x_{n}^{d}$ with $%
l\leq q$, $d<\upsilon ,$ $h\geq h_{1},$ and $s\in \lbrack 1,p]$ such that \ $%
h_{1}m_{0}+m_{s_{1}}+l_{1}m_{p}+d_{1}m_{n}=hm_{0}+m_{s}+lm_{p}+dm_{n},$ with
either $l\leq q^{\prime }$ or \ $d<\upsilon -w.$
\end{lemma}

\begin{proof}
For the case $q_{z}=0$ and $\varepsilon =0$ we have $q^{\prime }=q$ as $%
q^{\prime }=q-q_{z}-\varepsilon $, therefore the proof is done by the
previous lemma. Thus, we only need to work the proof for the case $q_{z}>0$
or $\varepsilon >0$. By the previous lemma the proof is reduced to $\sigma
=x_{0}^{h_{1}}x_{s_{1}}x_{p}^{l_{1}}x_{n}^{d_{1}}$ with $l_{1}\leq q$ and $%
d_{1}<\upsilon .$ Consider the following algorithm: \ \ 

\ \ \ \newline
Input $\sigma _{1}=x_{0}^{h_{1}}x_{s_{1}}x_{p}^{l_{1}}x_{n}^{d_{1}}$, $%
l_{1}\leq q$, $d_{1}<\upsilon .$\newline
Output $\sigma =x_{0}^{h}x_{s}x_{p}^{l}x_{n}^{d}$ with $l\leq q^{\prime }$
or $d<\upsilon -w$, and $s\in \lbrack 1,p]$.\newline
Set $i=1$.\newline
\textbf{WHILE} $l_{i}>q^{\prime }$ and $d_{i}\geq \upsilon -w$ \textbf{REPEAT%
}:$\newline
$\textbf{Step 1}: Write $d_{i}=a_{i}\left( \upsilon -w\right) +b_{i}$ with $%
a_{i}\in \mathbb{N}$ and $b_{i}\in \lbrack 0,\upsilon -w-1]$ (note $a_{i}>0$
since $d_{i}\geq \upsilon -w$). \newline
Write $l_{i}=a_{i}^{\prime }(q^{\prime }+1)+b_{i}^{\prime }$ with $%
a_{i}^{\prime }\in \mathbb{N}$ and $b_{i}^{\prime }\in \lbrack 0,q^{\prime
}] $ (note $a_{i}^{\prime }>0$ since $l_{i}>q^{\prime }$).\newline
Let $k_{i}=\min \{a_{i},a_{i}^{\prime }\}$. Note $k_{i}>0$ since $a_{i}>0$
and $a_{i}^{\prime }>0$. \ \ \newline
Reduce $\sigma _{i}$ with respect to $\psi _{p-r^{\prime }}$ $\ k_{i}$ times:%
\begin{equation*}
\sigma ^{\prime }=x_{0}^{h_{i}+k_{i}(\nu -1)}x_{s_{i}}x_{p-r^{\prime
}}^{k_{i}}x_{p}^{l_{i}-k_{i}(q^{\prime }+1)}x_{n}^{d_{i}-k_{i}\left(
\upsilon -w\right) }.
\end{equation*}%
\textbf{Step 2}: Reduce $\sigma ^{\prime }$ with respect to various $\xi
_{a,b}$ we get:%
\begin{equation*}
\sigma _{i+1}=x_{0}^{h_{i+1}}x_{s_{i+1}}x_{p}^{l_{i+1}}x_{n}^{d_{i+1}}.
\end{equation*}%
with some $l_{i+1}\leq l_{i}-k_{i}q^{\prime }\leq l_{i}$ (since $k_{i}>0$), $%
d_{i+1}=d_{i}-k_{i}\left( \upsilon -w\right) <d_{i}$, $h_{i+1}\geq h_{i}$,
and some $s_{i+1}\in \lbrack 1,p]$.\newline
Let $i=i+1$.\newline
\textbf{END\ LOOP}

\ \ \ 

It is clear that the sequence $l_{1},l_{2},\ldots $ is decreasing and the
sequence $d_{1},d_{2},\ldots $ is strictly decreasing. Hence the above
algorithm must terminate.
\end{proof}

\begin{lemma}
\label{Reduce Xs}Assume $W\neq \phi $ and $q^{\prime }=0$. Let $\sigma
=x_{0}^{h_{1}}x_{s_{1}}x_{n}^{d_{1}}$ with $d_{1}<\upsilon $, $h_{1}\geq 0,$
and $s_{1}\in \lbrack 1,p].$ Then $\sigma $ can be reduced to $%
x_{0}^{h}x_{s}x_{n}^{d}$ with\ $\left( s,d\right) \in V-W$ and $h\geq h_{1}$.
\end{lemma}

\begin{proof}
As $q^{\prime }=0$ then $\psi _{j}=\underline{x_{r^{\prime
}+j}x_{n}^{\upsilon -w}}-x_{0}^{\nu -1}x_{j}$. Also, either $q_{z}>0$ or $%
\varepsilon >0$ since $q>0$ and $q^{\prime }=q-q_{z}-\varepsilon $, so that $%
J=[0,p-r^{\prime }]$. Let $a=min\{t\geq 1\mid s_{1}-tr^{\prime }<r^{\prime
}\}$ and $b=min\{t\geq 1\mid d_{1}-t(\upsilon -w)<\upsilon -w\}$. Let $%
k=min\{a,b\}$. Reduce $\sigma $ with respect to $\psi _{s_{1}-ir^{\prime }}$
for $i=1,\ldots ,k$ we get $x_{0}^{h+k(\nu -1)}x_{s_{1}-kr^{\prime
}}x_{n}^{d_{1}-k(\upsilon -w)}$, with either $s_{1}-kr^{\prime }<r^{\prime }$
( hence $s_{1}-kr^{\prime }<u-z$) or $d_{1}-k(\upsilon -w)<\upsilon -w.$
Therefore, $(s_{1}-kr^{\prime },d_{1}-k(\upsilon -w))\in V-W$ as $s_{1}\leq
p<u$ and $d_{1}-k(\upsilon -w)<\upsilon $ for $k\geq 1$ since $%
d_{1}<\upsilon .$
\end{proof}

\begin{proposition}
\label{MainProp}Let $e_{0},\ldots ,e_{p},d_{1}\in \mathbb{N}_{0}$ with $%
e_{1}+\cdots +e_{p}\geq 1$. \textit{The monomial }$\alpha
=x_{0}^{e_{0}}\cdots x_{p}^{e_{p}}x_{n}^{d_{1}}$\textit{\ can }be reduced to 
$x_{0}^{h}x_{s}x_{p}^{l}x_{n}^{d}$\textit{\ (modulo }$\Phi $) \textit{with
respect to the grevlex monomial order,} \textit{such that }$%
\sum\limits_{i=0}^{p}e_{i}m_{i}+d_{1}m_{n}=hm_{0}+m_{s}+lm_{p}+dm_{n}$%
\textit{\ }with $(lp+s,d)\in V-W$.
\end{proposition}

\begin{proof}
We consider two cases $W=\phi $\ and $W\neq \phi $\ separately:

\ \ 

\underline{$W=\phi $:} $V-W$ $=V$\ . By Lemma~\ref{Lemma1} the proof is
reduced to the case $\sigma =x_{0}^{h}x_{s}x_{p}^{l}x_{n}^{d}$, with $l\leq
q $, $d<\upsilon $, and $s\in \lbrack 1,p]$. If $lp+s<u$ then $(lp+s,d)\in V$%
\ (since $d<\upsilon $), hence done. Otherwise, if $lp+s\geq u=qp+r$ then $%
l=q$ and $s\geq r$. Reduce $\sigma $ using $\varphi _{s-r}$ to $\sigma
_{2}=x_{0}^{h+\lambda -1}x_{s-r}x_{n}^{d+w}$. If $d+w<\upsilon $ then done
since $(s-r,d+w)\in V$ (since $s-r<p<u$). Otherwise, if $d+w\geq \upsilon $
then $w>0$, thus $z=0$ as $W=\phi $, hence $q_{z}=-1$ and $r_{z}=p$, thus $%
\theta =x_{n}^{\upsilon }-x_{0}^{\mu }$. Reduce $\sigma _{2}$ using $\theta $
to $x_{0}^{h+\lambda -1+\mu }x_{s-r}x_{n}^{d+w-\upsilon }.$ Observe $%
d+w-\upsilon <\upsilon $ since $d<\upsilon $ and $w<\upsilon $, hence $%
(s-r,d+w-\upsilon )\in V$ (since $s-r<p<u$), hence done.

\ 

\underline{$W\neq \phi $}: by Lemma~\ref{Lemma2} the proof is reduced to the
case $\sigma =x_{0}^{h}x_{s}x_{p}^{l}x_{n}^{d}$, with $l\leq q$, $d<\upsilon
,$ $s\in \lbrack 1,p]$, and with either $l\leq q^{\prime }$ or \ $d<\upsilon
-w$. Consider the following cases:\newline
\ \ \ \newline
\textbf{(1)} $lp+s<u-z$: then $(lp+s,d)\in V-W$ (since $d<\upsilon $), hence
done.\newline
\textbf{(2)} $lp+s\geq u-z$, consider three subcases:\ \ \ \newline
\qquad \textbf{(2a)} $d<\upsilon -w$ and $lp+s<u$: then $(lp+s,d)\in V-W$,
hence done.\newline
\qquad \textbf{(2b)} $d<\upsilon -w$ and $lp+s\geq u$: as $u=qp+r$ and $%
l\leq q$ then $l=q$ and $s\geq r$. Reduce $\sigma
=x_{0}^{h}x_{s}x_{p}^{l}x_{n}^{d}$ using $\varphi _{s-r}$ to $\sigma
^{\prime }=x_{0}^{h+\lambda -1}x_{s-r}x_{n}^{d+w}$. Note $(s-r,d+w)\in V$
because $s-r<p<u$ and $d+w<\upsilon .$\ If $(s-r,d+w)\in V-W$ then done.
Otherwise, if$\ (s-r,d+w)\notin V-W$ \ then necessarily $s-r\geq
u-z=q^{\prime }p+r^{\prime }$, therefore, $q^{\prime }=0$ because $s-r<p$.
Hence, done by Lemma~\ref{Reduce Xs}.\newline
\qquad \textbf{(2c)} $d\geq \upsilon -w$ : then by application of Lemma~\ref%
{Lemma2}, necessarily $l\leq q^{\prime }.$ Then $q^{\prime }p+r^{\prime
}=u-z\leq lp+s\leq q^{\prime }p+s$, so that $s\geq r^{\prime }$and $%
l=q^{\prime }.$ If we know that $s-r^{\prime }\in J$, we could use $\psi
_{s-r^{\prime }}$ to reduce $\sigma $. According to that we consider two
subcases:\ \ \newline
\qquad \textbf{(2c-1)} $q_{z}=0$ and $\varepsilon =0$: this implies $%
q=q^{\prime }$, therefore $q^{\prime }>0$. Under the assumptions, we have $%
l=q^{\prime }=q.$ Consider two cases:\newline
\underline{Case $s\geq r$:} reducing $\sigma $ with respect to $\varphi _{s-r%
\text{ \ }}$we get $\sigma ^{\prime }=x_{0}^{h+\lambda
-1}x_{s-r}x_{n}^{d+w}. $ As $d+w\geq \upsilon $ reduce $\sigma ^{\prime }$
with respect to $\theta $ we get $\sigma ^{\prime \prime }=x_{0}^{h+\lambda
+\mu -1}x_{s-r}x_{r_{z}}x_{n}^{d+w-\upsilon }.$ Note $s-r+r_{z}<p$ as $%
\varepsilon =0$ ($r>r_{z}$). Hence, using $\xi _{s-r,r_{z}}$ (or $\xi
_{r_{z},s-r}$ if $r_{z}\leq s-r$) $\sigma ^{\prime \prime }$ reduces to $%
x_{0}^{h+\lambda +\mu }x_{s-r+r_{z}}x_{n}^{d+w-\upsilon }$ with $%
(s-r+r_{z},d+w-\upsilon )\in V-W$ since $s-r+r_{z}<p<u-z=q^{\prime }p+r$ as $%
q^{\prime }>0$ and $d+w-\upsilon <\upsilon $. Hence, done.\newline
\underline{Case $s<r$:} in this case $0\leq s-r^{\prime }<r-r^{\prime
}=r_{z} $. Thus, $s-r^{\prime }\in J.$ Reduce $\sigma $ with respect to $%
\psi _{s-r^{\prime }\text{ \ }}$we get $x_{0}^{h+\nu -1}x_{s-r^{\prime
}}x_{n}^{d-(\upsilon -w)}$ with $s-r^{\prime }<p<u-z=q^{\prime }p+r$ as $%
q^{\prime }>0$. Hence, $\left( s-r^{\prime },d-(\upsilon -w)\right) \in V-W$
as $d-(\upsilon -w)<\upsilon $. Hence, done.\ \ \ \ \ \ \ \ \newline
\qquad \textbf{(2c-2)} $q_{z}>0$ or $\varepsilon >0$: hence $%
J=[0,p-r^{\prime }]$. If $q^{\prime }=0$ then done by Lemma~\ref{Reduce Xs}.
If $q^{\prime }>0$ then reduce $\sigma $ using $\psi _{s-r^{\prime }}$ to $%
x_{0}^{h+\nu -1}x_{s-r^{\prime }}x_{n}^{d-\left( \upsilon -w\right) }$. Note 
$\left( s-r^{\prime },d-\left( \upsilon -w\right) \right) \in V-W$ as $%
d-\left( \upsilon -w\right) <\upsilon $ and $s-r^{\prime }<p<u-z=q^{\prime
}p+r^{\prime }$ since $q^{\prime }>0$. Hence, done.
\end{proof}

\begin{proposition}
\label{SecondMainProp}\textit{Let }$\alpha =x_{0}^{e_{0}}\cdots
x_{p}^{e_{p}}x_{n}^{h}$\textit{\ and }$\beta =x_{0}^{d_{0}}\cdots
x_{p}^{d_{p}}x_{n}^{k}$\textit{\ with }$e_{0}$, $\ldots $, $e_{p}$, $h$, $%
d_{0}$, $\ldots $, $d_{p},$ $k\in \mathbb{N}_{0}$, with $e_{1}+\ldots
+e_{p}\geq 1$ and \ \ $d_{1}+\ldots +d_{p}\geq 1$ \textit{such that }$%
\sum\limits_{i=0}^{p}e_{i}m_{i}+hm_{n}=\sum%
\limits_{i=0}^{p}d_{i}m_{i}+km_{n} $. \textit{Then }$\alpha -\beta $\textit{%
\ can be reduced to zero modulo~}$\Phi $\textit{.}\ \ \ \ \ \newline
\end{proposition}

\begin{proof}
By Proposition~\ref{MainProp} $\alpha $\ can be reduced to\ $%
x_{0}^{c^{\prime }}x_{s^{\prime }}x_{p}^{l^{\prime }}x_{n}^{h+d^{\prime }}$%
\textit{\ }with $c^{\prime }\geq 0,(l^{\prime }p+s^{\prime },d^{\prime })\in
V-W$, and $s^{\prime }\in \lbrack 1,p]$. Furthermore, working with the same
procedures as in the proof of Proposition~\ref{MainProp} the monomial $%
x_{0}^{c^{\prime }}x_{s^{\prime }}x_{p}^{l^{\prime }}x_{n}^{h+d^{\prime }}$%
\textit{\ }can be reduced to $x_{0}^{c}x_{s}x_{p}^{l}x_{n}^{d}$\textit{\ }%
with $c\geq 0,(lp+s,d)\in V-W$, and $s\in \lbrack 1,p]$. Note that we have $%
cm_{0}+lm_{p}+m_{s}+dm_{n}=\sum\limits_{i=0}^{p}e_{i}m_{i}+hm_{n}=\sum%
\limits_{i=0}^{p}d_{i}m_{i}+km_{n}$. Hence, the proof follows by Lemma~\ref%
{Parameters} part (b).\ 
\end{proof}

\ \ \ 

\begin{proof}
(\textbf{of the main theorem, Theorem ~\ref{MainThm}}) We apply Buchberger
algorithm to show that all the S-polynomials, $S(f,g)$ with $f,g\in \Phi $,
reduces to $0$ with respect to $\Phi $ $\left( \text{modulo }\Phi \right) $.
Note that, by Lemma~\ref{Parameters}, all the binomials in the set $\Phi $\
are homogeneous with respect to the grading $wt(x_{i})=m_{i}$. Thus all the
S-polynomials are homogeneous. Also, it is easy to see that each such
S-polynomial has the form $\alpha -\beta $\ as in Proposition~\ref%
{SecondMainProp}\textit{. }This shows that $\Phi $ is a Groebner basis for
the ideal $inP.$

To complete the proof of Theorem~\ref{MainThm} we need to prove the
minimality of $\Phi $. The proof is as follows: it is clear that $LM(\theta
)\notin \langle LM(G-\{\theta \})\rangle $. By Lemma~\ref{Parameters} $%
w<\upsilon $ and by Remark~\ref{q>0-and-u>p} $q>0$, hence it is clear that $%
LM(\varphi _{i})\notin \langle LM(G-\{\varphi _{i}\})\rangle $ and $LM(\xi
_{i,j})\notin \langle LM(G-\{\xi _{i,j}\})\rangle $. If $q_{z}>0$ or $%
\varepsilon >0$ then $LM(\psi _{j})\notin \langle LM(G-\{\psi _{j}\})\rangle 
$ since $q^{\prime }<q$ because $q^{\prime }=q-q_{z}-\varepsilon $. If $%
q_{z}=0$ and $\varepsilon =0$ then $\{LM(\psi _{j})\mid 0\ \leq \ j\leq
r_{z}-1\}=\{x_{j}x_{p}^{q}x_{n}^{\upsilon -w}\mid r-r_{z}\ \leq \ j\leq
r-1\} $, hence $LM(\psi _{j})\notin \langle LM(G-\{\psi _{j}\})\rangle $.
\end{proof}

\ \ \ \ \ 

Patil and Singh (1990) constructed a generating set (but not minimal) for
the defining ideal $P$. We call the elements of this set \textquotedblleft
Patil-Singh generators\textquotedblright . This set of generators consists
of the same binomials in the set of Patil generators but with adding a few
more generators. The Patil-Singh generators are as follows: let $%
q,r,q_{z},r_{z},q^{\prime },r^{\prime },$ and $\varepsilon $ and be as
before,

\ \ \ 

\begin{tabular}{lll}
$\varphi _{i}$ & $=\underline{x_{i+r}x_{p}^{q}}-x_{0}^{\lambda
-1}x_{i}x_{n}^{w}$, & for $\ \ 0\leq i\leq p-r$; \ \ \ \  \\ 
$\psi _{j}$ & $=\underline{x_{r^{\prime }+j}x_{p}^{q^{\prime
}}x_{n}^{\upsilon -w}}-x_{0}^{\nu -1}x_{j}$, & for $\ \ 0\leq j\leq
p-r^{\prime }$; \\ 
$\theta $ & $=\underline{x_{n}^{\upsilon }}-x_{0}^{\mu
}x_{r_{z}}x_{p}^{q_{z}}$, &  \\ 
$\xi _{i,j}$ & $=\left\{ 
\begin{tabular}{ll}
$\underline{x_{i}x_{j}}-x_{0}x_{i+j}$, & $\text{if\ \ \ \ }i+j\leq p$; \\ 
$\underline{x_{i}x_{j}}-x_{i+j-p}x_{p}$, & $\text{if \ \ \ }i+j>p$,%
\end{tabular}%
\right. $ & for $\ \ 1\leq i\leq j\leq p-1$.%
\end{tabular}

\ \ \ \ 

The following theorem follows directly from Theorem~\ref{MainThm}.

\begin{theorem}
\label{SecondMainThm}\textit{The set }$\{\varphi _{i}\mid 0\leq i\leq
p-r\}\cup \{\theta \}$\textit{\ }$\cup $\textit{\ }$\{\xi _{i,j}\mid 1\leq
i\leq j\leq p-1\}$\textit{\ }$\cup $\textit{\ }$\{\psi _{j}\mid 0\ \leq \
j\leq p-r^{\prime }\}$\textit{\ is a Groebner basis for the ideal }$P$%
\textit{\ with respect to the grevlex monomial order with }$%
x_{0}<x_{1}<\cdots <x_{n}$\textit{\ with the grading }$wt(x_{i})=m_{i}$%
\textit{.}\ \ \ \ \ \ \ \ 
\end{theorem}

Note that Theorem~\ref{SecondMainThm} gives a Groebner basis with an easier
description but it gives up the minimality since if $q_{z}=0$ and $%
\varepsilon =0,$ then $r_{z}<r_{z}+p-r=p-r^{\prime }$and $\{LM(\psi
_{j})=x_{j+r^{\prime }}x_{p}^{q}x_{n}^{\upsilon -w}$ $\mid r_{z}\ \leq $\ $%
j\leq p-r^{\prime }\}=\{x_{j}x_{p}^{q}x_{n}^{\upsilon -w}$ $\mid r\ \leq $\ $%
j\leq p\}=x_{n}^{\upsilon -w}\{LM(\varphi _{i})=x_{j}x_{p}^{q}$ $\mid 0\
\leq $\ $j\leq p-r\}$.

\ \ \ 

Finally, we finish this article by noting that Patil-Singh generators do not
form a Groebner basis in all cases if we consider the grevlex monomial order
with the same grading as before but with\ $x_{0}>x_{1}>\cdots >x_{n}$ ( in
this case \textit{\ }$\prod\limits_{i=0}^{n}x_{i}^{a_{i}}>_{grevlex}\prod%
\limits_{i=0}^{n}x_{i}^{b_{i}}$\ if in the ordered tuple $%
(a_{1}-b_{1},\ldots ,a_{n}-b_{n})$\ the right-most nonzero entry is
negative). In the following we prove this and give an example.

\begin{remark}
\label{P-SnotGB}Assume $r<r_{z}<p$ (hence $\varepsilon =0$), $\lambda >1$,
and $w>0$. Then \textit{Patil-Singh generators are not a Groebner basis with
respect to the grevlex monomial ordering with }$x_{0}>x_{1}>\cdots >x_{n}$ 
\textit{and \ with the grading }$wt(x_{i})=m_{i}$\textit{\ .}
\end{remark}

\begin{proof}
First note the leading terms of the generators: $LT(\varphi
_{i})=x_{i+r}x_{p}^{q}$ if $w>0$ and $LT(\varphi _{i})=x_{0}^{\lambda
-1}x_{i}$\ if $w=0$, $LT(\psi _{j})=x_{0}^{\lambda +\mu -\varepsilon }x_{j}$%
, $LT(\theta )=x_{0}^{\mu }x_{r_{z}}x_{p}^{q_{z}}$, and $LT\left( \xi
_{i,j}\right) =x_{i}x_{j}$. Assume $r<r_{z}<p$ (hence $\varepsilon =0$), $%
\lambda >1$, and $w>0$. Consider $S(\xi _{1,r_{z}},\theta
)=x_{1}x_{n}^{\upsilon }-x_{0}^{\mu +1}x_{r_{z}+1}x_{p}^{q_{z}}$. Under the
above assumptions it is clear that none of the terms of $S(\xi
_{1,r_{z}},\theta )$ is a multiple of any of the leading terms of the
Patil-Singh generators.
\end{proof}

\begin{example}
Let $m_{0}=20,m_{1}=21,m_{2}=22,m_{3}=23,m_{4}=24,$ and $m_{5}=29$. Hence $%
p=4$. Let $P$\ be the kernel of the $K$-algebra homomorphism $\eta
:K[x_{0},\ldots ,x_{5}]\rightarrow K[t]$, defined by $\eta (x_{i})=t^{m_{i}}$%
. Recall the parameters in Lemma~\ref{Parameters}. It is easy to check that $%
\upsilon =3$, hence by the uniqueness condition we must have $\mu =2$, $%
q_{z}=1$, and $r_{z}=3$, thus $z=7$. For $1\leq i\leq 3$ note that in order
for $am_{4}+m_{i}-m_{0}$ to be in $\ \Gamma $ we must have $a\geq 2$. Note $%
g_{2p+1}=2(24)+21=2(20)+29$. Therefore, we conclude that $q=2$ and $r=1$,
thus $u=9$. Hence, $\lambda =2$, $w=1$, $r^{\prime }=2$, and $q^{\prime }=1$%
. Therefore, Patil-Singh generators are as follows: $G=\{\varphi _{i}\mid
0\leq i\leq 3\}\cup \{\psi _{j}\mid 0\leq j\leq 2\}\cup \{\theta \}\cup
\{\xi _{i,j}\mid 1\leq i\leq j\leq 3\}$ where $\varphi _{i}=\underline{%
x_{i+1}x_{4}^{2}}-x_{0}x_{i}x_{5}$, and $\psi _{j}=x_{j+2}x_{5}^{2}-%
\underline{x_{0}^{3}x_{j}}$, and $\theta =x_{5}^{3}-\underline{%
x_{0}^{2}x_{3}x_{4}}$ and $\xi _{i,j}$ as defined before with $p=4$. The set 
$G$ is not Groebner basis with respect to the grevlex monomial order with $%
x_{0}>x_{1}>\cdots >x_{5}$ and with the grading $wt(x_{i})=m_{i}$: consider $%
S(\theta ,\xi _{1,3})=x_{1}x_{5}^{3}-x_{0}^{3}x_{4}^{2}$. Note that neither
term of $S(\theta ,\xi _{1,3}))$ is a multiple of any of the leading terms
above.\vspace{0in}
\end{example}

\section*{Acknowledgement}

The author thanks Prof. Swanson I. for the useful discussions and comments
during the course of this work.

\ \ \ 

{\small Department of Mathematics and Statistics}

{\small Jordan University of Science and Technology}

{\small P O Box 3030, Irbid 22110, Jordan}

{\small Email address: iayyoub@just.edu.jo}

{\small \ \ \ \ \ \ \ \ }

\end{document}